\documentclass{gtart}
\usepackage{amsmath,amssymb,amscd}

\input gtmonout
\volumenumber{2}
\volumeyear{1999}
\volumename{Proceedings of the Kirbyfest}
\pagenumbers{291}{298}
\papernumber{15}
\received{19 January 1999}\revised{8 February 1999}
\published{20 November 1999}

\newtheorem{thm}{Theorem}
\newtheorem{prop}[thm]{Proposition}

\newtheorem{cor}[thm]{Corollary}

\theoremstyle{remark}
\newtheorem{rem}{Remark}

\theoremstyle{definition}

\title{On regularly fibered complex surfaces}
\author{D Kotschick}
\address{Mathematisches Institut, Universit\"at M\"unchen\\
Theresienstr.~39, 80333 M\"unchen, Germany}
\asciiaddress{Mathematisches Institut, Universitat Munchen,
Theresienstr. 39, 80333 Munchen, Germany}
\email{dieter@member.ams.org}

\begin{abstract}
We show that a compact complex surface which fibers
smoothly over a curve of genus $\geq 2$ with fibers of genus $\geq 2$
fibers holomorphically. We deduce an improvement of a result in~\cite{map},
and a characterisation of fibered surfaces with zero signature.
\end{abstract}
\asciiabstract{We show that a compact complex surface which fibers
smoothly over a curve of genus >=2 with fibers of genus >=
fibers holomorphically. We deduce an improvement of a result in
[D Kotschick, Math. Research Letters, {5} (1998) 227-234],
and a characterisation of fibered surfaces with zero signature.}


\begin{document}

\makeshorttitle

\centerline{\small\it Dedicated to Robion C~Kirby on the occasion of his
60$^{\text th}$ birthday.}

\section{Introduction}

In this paper we begin a study of complex structures on the 
total spaces $X$ of fiber bundles whose base $B$ and fiber
$F$ are compact orientable two--manifolds. We shall assume 
throughout that the genera $g(B)=g$ and $g(F)=h$ are at least
$2$. For fixed $g$ and $h$ we have infinitely many homotopy types
of orientable total spaces $X$, corresponding to conjugacy classes of
representations of the fundamental group of $B$ in the mapping
class group of $F$. By the Thurston construction all the total
spaces are symplectic, compare~\cite{map}.

Suppose now that $X$ admits some complex structure. The assumption
$g$, $h\geq 2$ implies in particular that $X$ is minimal and of 
general type. It has fixed topological Euler characteristic 
$c_2(X)=(2g-2)\cdot (2h-2)$ and therefore there are only finitely
many possible values for $c_1^2(X)$. By the boundedness results
of Moishezon and Gieseker
for the moduli space of surfaces of general type, we conclude that
among the infinitely many total spaces of surface bundles with
fixed $g$ and $h$ there are at most finitely many which admit a 
complex structure. In this paper we take the first step
or two towards characterising them.

In section~\ref{s:fiber} we shall prove a result which
implies that every complex structure as above admits a
regular holomorphic map $f$ to $B$ endowed with a suitable
complex structure. In fact, $f$ is in the same homotopy class
as the bundle projection. (This then makes the above finiteness 
result a consequence of Parshin--Arakelov finiteness, without
the need to appeal to Moishezon and Gieseker.) 

In section~\ref{s:app} we give some applications of the fibration 
criterion. We characterise the total spaces of surface bundles 
admitting complex structures under the additional assumption that 
the signature vanishes, and we also sharpen a result 
from~\cite{map}. This result concerns the minimal genus of a 
holomorphic representative for a given second cohomology class in 
the classifying space of the mapping class group of $F$, and thus 
bears on the complex analytic version of Problem~2.18 in Kirby's 
problem list~\cite{kirby}.

{\bf Acknowledgement}\qua I am grateful to F~Catanese and to D~Toledo for useful
comments.

\section{A fibration criterion}\label{s:fiber}

We shall say that a compact complex surface is {\it regularly
fibered} if it admits an everywhere regular holomorphic map
onto a smooth curve. It follows that the base curve and all the
fibers are compact, and that the surface is a smooth fiber 
bundle, though not usually a complex analytic bundle,
as the complex structure of the curves can vary.
We shall assume throughout that the genera of the base and
of the fiber are at least $2$.

Nontrivial examples of regularly fibered surfaces were first 
exhibited by Kodaira~\cite{kodaira} and by Atiyah~\cite{atiyah}, 
and later by many others, see eg~\cite{kas,hirz}.
Sometimes the terms Kodaira surface and, more often, Kodaira 
fibration are used to denote these surfaces. That terminology 
gives rise to confusion, not only because the term Kodaira
surface is more commonly used for certain non-K\"ahler complex
surfaces of Kodaira dimension zero, but also because some authors
seem to use Kodaira fibration to denote any regularly fibered
surface (with base and fiber of genus at least $2$) whereas others
implicitly restrict the term to mean only the examples of 
Kodaira~\cite{kodaira} (and maybe Kas~\cite{kas}) constructed as 
branched covers of products.

Here is the fibration criterion saying that any surface satisfying 
the obvious necessary conditions is indeed regularly fibered.

\begin{prop}\label{p:fiber}
Let $S$ be a compact complex surface whose fundamental group
fits into an extension
\begin{equation}\label{eq:ext}
1 \rightarrow \pi_1(F) \rightarrow \pi_1 (S) 
\stackrel{\pi}{\rightarrow} \pi_1(B)
\rightarrow 1 \ ,
\end{equation}
where $F$ and $B$ are closed oriented $2$--manifolds with genera
$g(F)=h\geq 2$ and $g(B)=g\geq 2$.
\begin{enumerate}
\item The topological Euler characteristic $e(S)\geq e(F)\cdot e(B)>0$.
\item The following conditions are equivalent:
\begin{itemize}
       \item $e(S)=e(F) \cdot e(B)$,
       \item $\pi$ is induced by a regular fibration of $S$ 
over $B$ endowed with a suitable complex structure,
       \item $S$ is aspherical.
\end{itemize}
\end{enumerate}
\end{prop}
\begin{proof}
Corresponding to the extension~\eqref{eq:ext} there is a 
fiber bundle $F \rightarrow X \rightarrow B$, with $X$ a 
classifying space for $\pi_1(S)$. As $X$ realises the smallest
possible Euler characteristic among all orientable $4$--manifolds 
with fundamental group $\pi_1(S)$, cf~\cite{groups}, we obtain 
$e(S)\geq e(X)=e(F)\cdot e(B)$, as claimed, where the last 
equality follows from the multiplicativity of the Euler 
characteristic in fiber bundles.

For the characterisation of the case of equality, notice that
if $S$ is regularly fibered, then it is aspherical by the 
homotopy exact sequence of the fibration. Further, if $S$ is
aspherical, then by the uniqueness of classifying spaces it is 
homotopy equivalent to $X$ and therefore has the same Euler
characteristic. Thus the crucial step is to show that the 
equality of Euler characteristics implies that $\pi$ is 
induced by a regular holomorphic map.

First of all, as $S$ is minimal with $c_2(S)=e(S)>0$ and with
$b_1(S)\geq 4$, it must be of general type. In particular, it is 
K\"ahler. By the theorem of Siu and Beauville, see Chapter~2 
of~\cite{book} and also~\cite{JY}, there is a surjective 
holomorphic map with connected fibers $f\co S\rightarrow C$, 
with $C$ a compact complex curve, such that the map $\pi$ 
in~\eqref{eq:ext} factors through $f_*$. This implies that 
$ker(f_* )$ is a (finitely generated) subgroup of $ker(\pi )=
\pi_1(F)$. Thus $ker(f_* )$ is the fundamental group of an 
orientable surface $\overline{F}$ which is a covering of $F$. If 
$\overline{F}$ were noncompact, $ker(f_* )$ would be a free group, 
contradicting the fact that $\pi_1(S)$ has cohomological dimension 
$4$. Thus $\overline{F}$ is compact, with 
\begin{equation}\label{eq:1}
g(\overline{F})\geq g(F)=h \ .
\end{equation}

On the other hand, denoting the generic fiber of $f$ by $F'$, we 
have that $ker(f_* )=\pi_1(\overline{F})$ is a quotient of $\pi_1(F')$ 
and so $g(\overline{F})\leq g(F')$. Now by the theorem of 
Zeuthen--Segre a singular fiber makes a positive contribution
to the Euler characteristic, so we have 
$e(S)\geq (2g(C)-2)\cdot (2g(F')-2)$, so that
$e(S)=(2g(B)-2)\cdot (2g(F)-2)$ and $g(C)\geq g(B)=g$ imply
\begin{equation}\label{eq:2}
g(\overline{F})\leq g(F')\leq g(F)=h \ .
\end{equation}

Combining~\eqref{eq:1} and~\eqref{eq:2}, we conclude that
$g(\overline{F})=g(F')=g(F)=h$ and therefore $g(C)=g(B)=g$.
Thus $C$ gives a complex structure on $B$ and $f$ is a 
holomorphic map inducing $\pi$. As we are in the case of 
equality for the Zeuthen--Segre inequality, $f$ must be 
everywhere regular.
\end{proof}

In Chapter 2 of~\cite{book} the genus $g(M)$ of a compact K\"ahler 
manifold $M$ was defined. This is the maximal genus of a compact
curve $C$ onto whose fundamental group $\pi_1(M)$ surjects. The 
Siu--Beauville theorem shows that if $\pi\co\pi_1(M)\rightarrow 
\pi_1(C)$ is any surjective homomorphism with $g(M)=g(C)$, then 
$\pi$ is induced by a holomorphic map with connected fibers. This 
conclusion does not necessarily hold if $g(M)>g(C)$, although 
surjective homomorphisms will exist in abundance. An interesting 
aspect of the proof of Proposition~\ref{p:fiber} is that it shows
the homomorphism $\pi$ in~\eqref{eq:ext} is induced by a holomorphic 
map with connected fibers although the genus of $S$ may very well be 
larger than $g(B)=g$: just take a trivial extension with $h=g(F)>g(B)=g$.

\begin{rem}
A version of the second part of Proposition~\ref{p:fiber} has been 
proved idependently by Hillman~\cite{hillman}, but his proof is more 
complicated. He begins by using the work of Gromov and of
Arapura--Bressler--Ramachandran on $L^2$--cohomology (see
Chapter 4 of~\cite{book}) to produce a holomorphic map to a
curve. In~\cite{kapovich} an extension of the argument is proposed
in the case where the kernel of $\pi$ in~\eqref{eq:ext} is not 
assumed to be a surface group, but can be any finitely presentable 
group. It turns out that this more general result can be 
deduced from Proposition~\ref{p:fiber} or from the result 
of~\cite{hillman} using standard arguments on the cohomology
of Poincar\'e duality groups~\cite{hillman2}.
\end{rem}

As an immediate consequence of Proposition~\ref{p:fiber} we have:

\begin{cor}\label{c:1}
If $S$ is any compact complex surface homotopy equivalent to a 
surface bundle $X$ over a surface with base and fiber of genera at 
least $2$, then $S$ is regularly fibered and is diffeomorphic to $X$.
\end{cor}

This generalises results of Kas~\cite{kas} and of Jost--Yau~\cite{JY}
who showed that deformations of Kodaira's examples~\cite{kodaira}
are regularly fibered. In those examples one obtains a description
of a component of the moduli space in terms of moduli spaces of 
curves underlying the construction. In the general case,
the corollary says that all components of the moduli space of 
complex structures on this particular manifold are made up of 
regularly fibered surfaces, but there is no direct description in
terms of the moduli of curves.

\section{Applications}\label{s:app}

The original motivation for studying regularly fibered surfaces
was that they provide examples of smooth fibre bundles for which
the signature is not multiplicative~\cite{atiyah,kodaira}, in
this case that just means non-zero. In~\cite{map} we proved some
bounds on the signatures of surface bundles over
surfaces. We shall now slightly improve Theorem 3 of~\cite{map}:

\begin{thm}\label{t:MY}
Let $X$ be a surface bundle over a surface, with the genera of
the base and the fiber $\geq 2$. 
If $X$ admits a complex structure (not necessarily 
compatible with the orientation), or an Einstein metric, then
\begin{equation}\label{eq:MY}
3\vert\sigma (X)\vert < e(X) \ .
\end{equation}
\end{thm}
\begin{proof}
Suppose $X$ admits a complex structure. After possibly
reversing the orientation, we may assume that the complex 
structure is compatible with the orientation. 

The argument in~\cite{map} was as follows:
$X$ is a minimal surface of general type for which the 
underlying manifold endowed with the other, non-complex, 
orientation is symplectic and therefore has non-zero 
Seiberg--Witten invariants. Thus Theorem~1 of~\cite{georient} gives 
$$
\sigma (X)\geq 0 \ . 
$$
This, together with the Miyaoka--Yau inequality 
$$
3\sigma (X)\leq e(X) \ ,
$$ 
implies $3\vert\sigma (X)\vert\leq e(X)$.

We can now reach the same conclusion in a different way, and we
can also show that the inequality must be strict, as claimed.
By Corollary~\ref{c:1} the surface $X$ is regularly fibered, 
so that the non-negativity of the signature follows from
Arakelov's theorem. Moreover the Miyaoka--Yau inequality is 
strict for regularly fibered surfaces, as proved by Liu~\cite{liu}.

Suppose that $X$ admits an Einstein metric. As it is 
also symplectic, it has non-zero Seiberg--Witten
invariants and by the result of~\cite{LeBrun}
satisfies $3\sigma (X)\leq e(X)$. The same argument
for the manifold with the other orientation gives
$-3\sigma (X)\leq e(X)$. It remains to exclude the case of
equality. 

Suppose that $X$ admits an Einstein metric and that for
a suitable choice of orientation $3\sigma (X) = e(X)$.
Then Le~Brun~\cite{LeBrun} showed that the Einstein metric
must be K\"ahler--Einstein, so that $X$ must be a complex surface
which is a ball quotient. But then by the above argument for the 
complex case, we have a contradiction.
\end{proof}

We now return to the issue of characterising those 
surface bundles over surfaces which admit complex structures.
Here is such a characterisation in the easiest case,
when the signature of the total space vanishes.

\begin{thm}\label{t:zero}
Let $X$ be the total space of a surface bundle over a surface, 
with the genera of the base $B$ and of the fiber $F$ at least $2$.
The following are equivalent:
\begin{enumerate}
\item $X$ admits a complex structure and has zero signature,
\item the monodromy representation $\rho\co\pi_1(B)\rightarrow
\Gamma_h$ has finite image\footnote{$\Gamma_h$ is the mapping
class group of the fiber $F$, where $h=g(F)$.}.
\end{enumerate}
\end{thm}
\begin{proof}
Suppose $X$ admits a complex structure, then by Proposition~\ref{p:fiber}
$X$ is regularly fibered. If the signature vanishes, we are in the
borderline case of Arakelov's theorem, which says that the signature is
nonnegative, and is zero only if all the fibers are isomorphic,
so the fibration is isotrivial. In this case we can pull 
back the fibration to a finite cover of $B$ to obtain a product. This 
implies that the kernel of the monodromy representation has finite 
index in $\pi_1(B)$.

Conversely, assume that we have a bundle with finite monodromy.
Then it must have zero signature.
By the positive resolution of the Nielsen realisation 
problem~\cite{kerk} we can choose a complex structure on $F$
and a lift of the monodromy representation to the diffeomorphism
group of $F$ so that the monodromy acts by complex analytic
diffeomorphisms of $F$. Fixing an arbitrary complex structure
on $B$, we obtain a complex structure on $X$ by viewing it
as $(F\times \tilde{B})/\pi_1(B)$, where $\pi_1(B)$ acts on
$\tilde{B}$ by deck transformations and acts on $F$ through the 
chosen lift of the monodromy representation to the diffeomorphism 
group.
\end{proof}

\begin{rem}
Under the conditions of the theorem $X$ is finitely covered
by a product, and is uniformised by the polydisk, compare~Theorem 1
in~\cite{georient}. Thus Theorem~\ref{t:zero} is related to
Catanese's characterisation~\cite{cat} of complex surfaces 
finitely covered by products.
\end{rem}

By the finiteness results for the complex case, Theorem~\ref{t:zero}
has the following immediate consequence:

\begin{cor}
For fixed $g=g(B)$ and $h=g(F)$, both $\geq 2$, there are only
finitely many conjugacy classes of representations $\rho\co
\pi_1(B)\rightarrow\Gamma_h$ with finite image.
\end{cor}

This is a considerable strengthening of the following result
of Harvey~\cite{harvey}:
\begin{cor}
The finite subgroups of the mapping class group $\Gamma_h$ fall into
finitely many conjugacy classes.
\end{cor}
\begin{proof}
Every finitely generated subgroup of $\Gamma_h$ is the monodromy group
of a surface bundle of zero signature over some base $B$, where the
genus of $B$ can be taken to be the number of generators of the 
subgroup, see Proposition 4 of~\cite{map}. As the order of the finite 
subgroups of $\Gamma_h$ is bounded\footnote{by $84(h-1)$}, we have an
a priori bound on $g=g(B)$ and can apply the previous corollary.
\end{proof}
It is clear that the first corollary is much stronger than the 
second one, as there are usually many different monodromy 
representations with the same image, compare section 3 of~\cite{map}.

Finally, note that every surface bundle with fibers of genus
$2$ has zero signature, and so is covered by Theorem~\ref{t:zero}.
In the higher genus case there are always bundles of non-zero
signature, and for these a characterisation of the monodromy
representations arising from complex surfaces is not yet available.
We shall return to this in a future paper. Here we just remark
that this problem need not be the same as trying to decide
which extensions of surface groups by surface groups are 
K\"ahler groups. If a surface bundle admits a complex structure, 
then its fundamental group is K\"ahler and in fact projective.
However, it is possible that there are surface bundles which
admit no complex structure but still have K\"ahler fundamental
groups.

\Addresses\recd

\end{document}